\documentclass[11pt]{article}
\usepackage{geometry}                
\usepackage{graphicx}
\usepackage{amssymb}
\usepackage{amsmath}
\usepackage{epstopdf}
\usepackage{authblk}
\usepackage{tikz}
\usepackage{amsthm}

\usepackage[all]{xy}
\title{Finding meaningful and workable applied mathematics problems in science}
\author{Yue Wang}

\date{}                                           

\begin{document}
\maketitle

\begin{abstract}
In this short review, I will summarize my research experience in three fields in applied mathematics: mathematical biology, applied probability, and applied discrete mathematics. Specifically, I will show how each project was initiated, and what wrong approaches were applied. Such details are important in learning how to do research, but they cannot be read out from research papers. I wish that students and junior researchers in applied mathematics could learn a lesson from this summary.
\end{abstract}

\section{Introduction}

I have been working on applied mathematics since 2012. Although my research might not be very attractive, my research experience should be helpful to students and junior researchers in applied mathematics, in both positive (what to do) and negative (what to avoid) sense. The main focus of this review is to present the story behind each paper: how did I find this problem, and what right/wrong approaches did I apply? Such lessons cannot be learned by just reading papers. I would be happy to read this type of review if I were much younger. This is the motivation of writing this summary. In this review, I will summarize my research in mathematical biology (Section~\ref{s2}), applied probability (Section~\ref{s3}), and applied discrete mathematics (Section~\ref{s4}).

For the convenience of readers, I decide to move the final discussions here. The following are some lessons you might learn from this summary, and you can directly jump to the corresponding section: After proving some results, try to strengthen them by considering a more general or more realistic setting. A small improvement might be worth doing if this improvement is necessary in the real world (Subsection~\ref{2.1}). As a researcher in applied mathematics, sometimes it is inevitable to work in an unfamiliar field. In this case, it is important to find a collaborator in the corresponding field. Therefore, it is always good to have more friends with different backgrounds (Subsection~\ref{2.2}). After solving a problem, think in the opposite direction, so that your results can be used to formulate another problem (Subsection~\ref{2.3}). Before tackling a hard problem, search for related papers and talk to different experts to avoid reinventing the wheel (Subsection~\ref{3.1}). When writing the paper, think about your potential readers (also reviewers), especially if the paper is related to multiple fields (Subsection~\ref{3.1}). When reading papers that propose fancy methods, check them in weird situations and try to construct counterexamples (Subsections~\ref{3.2},~\ref{3.3}). After developing new methods, always try your best to apply your methods to some real data, or at least simulated data. Otherwise, this will always be questioned by reviewers (Subsection~\ref{2.2}). Besides, to evaluate or compare algorithms, the best way is to implement and test them, since your theoretical analyses might overlook some factors. After learning something outside of your fields, try to apply your mathematical thinking mode to analyze it (Subsection~\ref{4.1}). Improve your coding skill. You need to implement the methods you develop, which is sometimes nontrivial, especially if you need to design fast algorithms (Subsection~\ref{4.2}). Otherwise, you need a collaborator for this. I hope that you can find meaningful and workable problems, and avoid the mistakes I made.

\section{Mathematical biology}
\label{s2}
\subsection{Population dynamics}
\label{2.1}

In some types of cancer, there are multiple phenotypes that can convert into each other. It has been observed that starting from any one phenotype, after some time, the final proportions of different phenotypes always converge to the same constants. This phenomenon is called phenotypic equilibrium. I read a paper that uses a Markov chain model to explain this phenomenon \cite{gupta2011stochastic}. This model has some flaws in mathematics and biology. I built a linear ODE model and proved the phenotypic equilibrium phenomenon as the convergence to the unique attracting fixed point \cite{zhou2014multi}. This work also intrigued other works that study multiple phenotypes in cancer \cite{niu2015phenotypic,chen2016overshoot}. Later, I built a branching process model, which makes more sense in biology. I tried to prove a law of large numbers in this branching process model to explain the phenotypic equilibrium phenomenon. I found some results in a probability paper \cite{janson2004functional}, but those results need a certain condition that all phenotypes form a single communicating class, which might not hold in reality. Based on the results in this paper, I used some complicated arguments to remove this condition and proved the law of large numbers I wanted \cite{jiang2017phenotypic}.

I collaborated with a wet lab that cultivated leukemia cells. There were many cell populations, starting from 10/4/1 cells, and the cell area of each population was measured every day. Due to some technical issues, the cell area data were not quite accurate for small populations. Therefore, the analysis has to be done only for relatively large populations. I processed the data and drew different figures. Finally, I found the best figure that illustrates the most important feature. In this figure, I translated all growth curves, so that they start at the time point of reaching 2500 cells. In this figure, all populations from 10 initial cells grew fast, while some populations from 1 initial cell grew much slower. I built many toy models and found the most probable one, which has multiple phenotypes with different growth rates. Although from the same cell line, there might be multiple phenotypes \cite{qian2020counting}. The difference in growth rates can be simply explained: starting from 10 cells, it is very likely to have at least one fast cell, which will dominate, and the overall growth rate is high; starting from 1 cell, it might be of the slow type, and the overall growth rate is low \cite{wang2018some,wang2023multiple}. For researchers in mathematics, communication and collaboration with experimental biologists might not be smooth at the beginning, since the thinking modes are different. I also participated in other projects on population dynamics \cite{angelini2022model,wang2022modelling,dessalles2021naive,xia2023age}.

\subsection{Gene regulation}
\label{2.2}
Genes are transcribed to mRNAs, which are then translated to proteins. Some proteins can affect the expression of other genes (mutual regulation) or their own corresponding genes (autoregulation). Since it is extremely difficult to directly determine gene regulatory relations through biochemical methods, there have been many methods to infer gene regulation \cite{zhou2021dissecting,bocci2022splicejac}. I read some papers on gene regulation inference, and found that the settings are diverse: different methods use different models, which are based on different assumptions. Besides, different methods may require different data types. For this complicated situation, I gradually came up with a classification system. The data types can be classified by four dimensions: measuring gene expression or phenotype; measuring at single-cell level or bulk level; measuring at one time point or multiple time points; measuring with or without interventions. For single-cell data at multiple time points, we need to distinguish whether one can measure the joint probability distribution, or just marginal probability distributions. In total, there are 20 data types. For each data type, I summarized known methods, or developed novel methods, or proved that data type could not be used to infer gene regulation. I implemented those new methods, and found that one method is too sensitive to small perturbations, thus not applicable to real data. My original plan was just to develop a novel inference method, but finally ended with a very long paper \cite{wang2022inference}. To finish such a long paper, I collaborated with a friend in experimental biology, who wrote the biological part. My ambition was to make this paper useful to researchers with different backgrounds. Therefore, I managed to modularize different sections, and added a figure indicating which sections are dependent. Then mathematicians and biologists can skip different sections to reach what they are interested in.

When studying the inference for mutual regulation, I also found an idea to infer the existence of autoregulation. In a Markov chain model of gene expression, if there is no autoregulation, then for the expression level, the variance should be larger than the mean \cite{wang2022ar}. Thus if the mean is larger than the variance, there might be autoregulation. Since this model has infinitely many states, I failed to prove what I wanted. Therefore, I sought help from a friend in pure mathematics. We collaborated and finished all the proofs. As an applied probability paper, the results are not deep enough; as a mathematical biology paper, the inference methods are not implemented. A common and natural question from reviewers for a methodology paper is its application. Thus I had to implement my methods and applied them to some experimental data.

\subsection{Developmental biology}
\label{2.3}
Consider tissue transplantation experiments, where a piece of donor tissue is transplanted to a host tissue. During the development of the host, we observe the fate of the transplanted tissue. The results can be normal/abnormal, host-like/donor-like. I got some summarized experimental results, and the goal is to develop some methods to infer the unknown experimental results. Since the data amount is too small, I tried several fancy models and failed. Finally, inspired by the energy function in the Ising model, I adopted a straightforward approach: we can estimate the similarities between experiments. Design a simple penalty function. For each guess for the unknown results, calculate the total penalty: a pair of similar experiments with different results (whether known or from guess) gets a higher penalty. The guess with the minimal penalty is the most probable guess \cite{wang2021inference}. After inferring unknown results, I considered this problem from the opposite direction: if we want to use this method to infer experimental results, and we can choose what experiments to perform, what is the optimal design that minimizes the cost? This problem turns out to be a coloring problem in $n$ dimensional lattice. With the help of a friend in pure mathematics, I solved this experimental design problem.

I once took over an unfinished project on positional information in developmental biology. In a developing embryo, cells at different positions seem to know where they are and develop into different (yet organized) tissues. To explain this, biologists coined a term, positional information \cite{wolpert2016positional}. This notion is used by different biologists in different meanings. I analyzed carefully and proposed three criteria for defining what is positional information. I also proved that no criterion is redundant, since removing any one criterion allows some non-positional information being included. In this logical way, positional information is clearly defined. Finally, I finished a mathematically styled paper without any mathematics \cite{wang2020biological}. I also contributed to another project in developmental biology \cite{wang2020model}.

\section{Applied probability}
\label{s3}
\subsection{Thermodynamics in stochastic processes}
\label{3.1}

My Ph.D. group has a tradition of working on thermodynamics of stochastic processes \cite{ma2015universal,ye2016stochastic,qian2019ternary,qian2020kinematic,yang2021potentials,cheng2021stochastic,cheng2021asymptotic}. I worked on a problem of comparing the thermodynamic properties of a Markov chain and its covering space. This problem is related to topology and group theory, which I did not know very well. After trying many different wrong directions, I found the right theorem that I should prove. That theorem looks correct, but I got stuck for a long time, and tried to seek advice from different mathematicians. Finally I realized that I was terrified by this problem, and a very loose inequality is enough. 

Later, I planned to generalize this result to more general covering spaces. I got stuck, asked many mathematicians, took a related course. Finally, I accidentally learned that a previous paper \cite{kaimanovich2002boundary} has covered all what I wanted to do, and even some results that I had proved. Therefore, my plan of generalization failed disastrously. The authors of that paper gave a very different name to the object I studied, so that I did not find it when I was searching for references. Besides, this small field is mostly studied in Europe, so that the mathematicians I asked did not know this paper. Fortunately, their methods are different from mine, so that my paper finally got published \cite{wang2020mathematical}. The publication of this paper is also difficult. This paper is related to statistical physics, probability, and algebraic topology. A reviewer might be unfamiliar with one of these subjects and get lost. With my techniques in probability, I also proved theorems in finite Markov chains \cite{wang2022discrete}.

\subsection{Causal inference}
\label{3.2}
Once my advisor (who does not work on causal inference) attended an interesting talk on quantifying causal effects, and later sent me some related papers. Traditional causal quantities have problems when two possible causes are very similar. There are some new methods that use the slight difference to enhance the performance \cite{janzing2013quantifying,zhao2016part}. From their definitions, when two possible causes are identical, their quantities are not defined. I tried to use continuation, namely adding small perturbations on the joint distribution to make the quantities defined, and then decreasing the perturbations to zero. This opened a can of worms: near the distribution where those quantities are not defied, such quantities fluctuate violently. Therefore, when two possible causes are identical (or more generally, the Markov boundary is not unique), those quantities are essentially ill-defined. I tried to develop new causal quantities that avoid this problem but failed. Then I realized that the problem itself might be ill-posed. I proposed three criteria for a causal quantity, and proved that they are incompatible when the distribution has multiple Markov boundaries \cite{wang2020causal}. Since distributions with multiple Markov boundaries are problematic, I also designed fast algorithms to detect them from finite data. Before implementing these algorithms, I did some theoretical analyses and argued that one algorithm is the best. After implementing them, I found that another algorithm is the best, since I underestimated some factors in my theoretical analyses. Since neither my advisor nor myself has training in causal inference, I had to spend a long time to find a collaborator. In such highly mathematical fields, guidance from an expert is essential both for controlling the research direction and writing in the correct style.

\subsection{Reinforcement learning}
\label{3.3}
Once a friend in operations research sent me some papers in online learning, where a player repeatedly chooses an action and receives a stochastic reward (also partially learns the system parameters) The goal of the player is to apply a wise policy to maximize the cumulated reward. I tried to understand a fancy paper, but realized that the measurement to evaluate the performance of a policy is problematic. The measurement only considers the performance in the worst case, which leads to some strange behaviors. Specifically, concrete historical data cannot help in locating the optimal policy. This field is quite friendly to someone with the background in probability, and I just spent two months to work out this problem \cite{wang2021measuring}. Nevertheless, such negative results are quite hard to get published. Later, I also studied a situation where the historical data are polluted \cite{wang2023online}. The problem of this paper is that the model is too simple (so that the optimal policy can be explicitly calculated), but more complicated models can be only studied numerically, which does not provide enough insights.

\section{Applied discrete mathematics}
\label{s4}
\subsection{Inheritance law}
\label{4.1}
I have a paper in inheritance law \cite{wang2022impossibility}, although I never take a course in law or even social science. This project has an unexpected initiation. Once I read a story about heritage succession on social media: Alice's father passed away before her paternal grandfather. Thus a portion of her father's heritage went to her paternal grandfather, and then (after the death of her paternal grandfather) went to her paternal uncles and aunts. If her father passed away after her paternal grandparents, then her uncles and aunts would not inherit from her father. The reason is that in China, the Civil Code regulates that parents, children, and spouse can equally split the heritage.

From this story, I learned that in China's inheritance law, different orders of death might lead to different inheritance results. If we regard death and inheritance as an operator on the wealth of different persons, then this operator is non-commutative. This leads to a natural question: what if multiple relatives died in the same disaster, and the order of death is unknown? 

I asked a friend in legal studies, and learned that this situation has been considered in China's inheritance law. The order of death is determined by relative relations: if one dead person has no living inheritable relatives, this person is presumed to die first; if both dead persons have living inheritable relatives, and they are of the same generation, they are presumed to die simultaneously and do not inherit each other; otherwise, the one of the senior generation is presumed to die first. For example, sisters are presumed to die simultaneously, and son is presumed to die after father.

This stipulation of death order does not look very rigorous. First, the term ``generation'' might not be well-defined. For example, Nobel laureate Gabriel Garc\'ia M\'arquez married his uncle's sister-in-law. In this case, how should one compare the generations of Gabriel, his wife, and his uncle?

Even without the trouble of ``generation'', the stipulation of death order still does not look right. Using some insights in discrete mathematics, I constructed a counterexample, where for multiple dead persons, the order of death is partially known, and the stipulated order for the unknown ones are contradicted with the known orders. 

Although that counterexample is not natural and not appreciated by that friend in legal studies, I started to find a way to assign the order of death, which does not lead to contradictions. After some trials, I found no good solution, and started to suspect whether such a solution exists. If we cannot assign the order of death, there is another approach: the need for assigning the order of death is from the fact that the inheritance law is non-commutative. I tried to construct a commutative inheritance law and also failed. 

Since neither approach worked, I guessed that the order of death problem cannot be solved. I proposed some reasonable criteria for an inheritance law. Then I proved that under such criteria, no inheritance law is commutative, and any assignment of the order of death can lead to contradictions. The proofs just need some basic tricks in discrete mathematics. In sum, any inheritance law is incomplete. Proposing reasonable criteria and prove that they are incompatible is a standard approach to derive impossibility results, similar to Arrow's impossibility theorem in voting \cite{maskin2014arrow}. Certainly, in practice, the order is reversed: I first considered how to prove such impossibility results, then tried to propose necessary criteria and remove redundant criteria.

A commutative inheritance law is incompatible with the three criteria I proposed (gender equality, inheritance rights of children, and conditional inheritance rights of parents). Since assigning the order of death is impossible, I switched to construct a commutative inheritance law that only violates one criterion. I proved that there are only four commutative inheritance laws that only violate one criterion. Therefore, that criterion is the price for a well-defined inheritance law.

After proving this theorem, I started to write a paper. To fill in the introduction section, I checked the inheritance laws in different countries. In the French Civil Code, it is stipulated that if one cannot determine the order of death for two persons, they do not inherit each other. This means that the order of death is not really assigned (incomparable). After learning this, I realized an embarrassing fact: I only proved that assigning the order of death leads to contradictions, but this French approach does not introduce new order and leads to no contradiction. Thus it is a perfect solution to the order of death problem. Then I modified my proof to include this situation, and the final result becomes: the French approach is the only valid solution. The results that only four commutative inheritance laws only violate one criterion also become useless in practice. Since they are of some theoretical interests, I moved them to the appendix.

This paper is in an awkward situation: researchers in legal studies do not understand the mathematical proofs, and applied mathematicians are not interested in law. Therefore, it is extremely difficult to find a proper journal and proper reviewers (although finally succeeded). Also, it is hard to seek guidance or find a collaborator.

\subsection{Algorithm design}
\label{4.2}
Although not well trained in computer science, I used my techniques in discrete mathematics to design algorithms and prove the correctness. One project considers finding transposons: given some annotated gene sequences from different individuals, such as $\{1,2,3,4,5,6\}$, $\{1,2,5,3,4,6\}$, $\{1,2,3,4,6,5\}$, how should one find the genes that change their locations? In this example, obviously gene $5$ changes its location. However, we can also explain the difference by assuming genes $3,4,6$ change their locations. Certainly, we can minimize the number of operations needed to transform all sequences into the same. I came up with a more concise solution: determine the longest common subsequence, and the complement is all the genes that translocate. With the help from a friend in computer science, I implemented a fast algorithm, which was used to compare the genes of bacteria \cite{kang2014flexibility}. Although determining the longest common subsequence is a classic problem in computer science, I also considered the problem with multiple longest common subsequences, which made a new paper \cite{wang2023longest}.

Another algorithm design work needs to compare two rooted trees with possibly repeated labels, where there is no order for different children vertices. I read many related papers but found no method that can be directly applied to this situation. I first designed a ``best match'' metric, which calculates the minimal distance between all possible transformed trees. However, I thought this metric could not be calculated in polynomial time, and I switch to another metric, which is fast to calculate. Later, I found that the best match metric could be calculated in quadratic time. Therefore, I contained both metrics in my paper, although the best match metric is significantly better \cite{wang2022two}. These metrics were used to compare the early development of plants. 

\bibliographystyle{acm}
\bibliography{am}

\begin{thebibliography}{10}

\bibitem{angelini2022model}
{\sc Angelini, E., Wang, Y., Zhou, J.~X., Qian, H., and Huang, S.}
\newblock A model for the intrinsic limit of cancer therapy: {Duality} of
  treatment-induced cell death and treatment-induced stemness.
\newblock {\em PLOS Computational Biology 18}, 7 (2022), e1010319.

\bibitem{bocci2022splicejac}
{\sc Bocci, F., Zhou, P., and Nie, Q.}
\newblock {spliceJAC}: transition genes and state-specific gene regulation from
  single-cell transcriptome data.
\newblock {\em Molecular Systems Biology 18}, 11 (2022), e11176.

\bibitem{chen2016overshoot}
{\sc Chen, X., Wang, Y., Feng, T., Yi, M., Zhang, X., and Zhou, D.}
\newblock The overshoot and phenotypic equilibrium in characterizing cancer
  dynamics of reversible phenotypic plasticity.
\newblock {\em Journal of Theoretical Biology 390\/} (2016), 40--49.

\bibitem{cheng2021stochastic}
{\sc Cheng, Y.-C., and Qian, H.}
\newblock Stochastic limit-cycle oscillations of a nonlinear system under
  random perturbations.
\newblock {\em Journal of Statistical Physics 182}, 3 (2021), 1--33.

\bibitem{cheng2021asymptotic}
{\sc Cheng, Y.-C., Qian, H., and Zhu, Y.}
\newblock Asymptotic behavior of a sequence of conditional probability
  distributions and the canonical ensemble.
\newblock {\em Annales Henri Poincar{\'e} 22}, 5 (2021), 1561--1627.

\bibitem{dessalles2021naive}
{\sc Dessalles, R., Pan, Y., Xia, M., Maestrini, D., D'Orsogna, M.~R., and
  Chou, T.}
\newblock How naive {T}-cell clone counts are shaped by heterogeneous thymic
  output and homeostatic proliferation.
\newblock {\em Frontiers in immunology 12\/} (2021).

\bibitem{gupta2011stochastic}
{\sc Gupta, P.~B., Fillmore, C.~M., Jiang, G., Shapira, S.~D., Tao, K.,
  Kuperwasser, C., and Lander, E.~S.}
\newblock Stochastic state transitions give rise to phenotypic equilibrium in
  populations of cancer cells.
\newblock {\em Cell 146}, 4 (2011), 633--644.

\bibitem{janson2004functional}
{\sc Janson, S.}
\newblock Functional limit theorems for multitype branching processes and
  generalized {P{\'o}lya} urns.
\newblock {\em Stochastic Processes and their Applications 110}, 2 (2004),
  177--245.

\bibitem{janzing2013quantifying}
{\sc Janzing, D., Balduzzi, D., Grosse-Wentrup, M., and Sch{\"o}lkopf, B.}
\newblock Quantifying causal influences.
\newblock {\em The Annals of Statistics 41}, 5 (2013), 2324--2358.

\bibitem{jiang2017phenotypic}
{\sc Jiang, D.-Q., Wang, Y., and Zhou, D.}
\newblock Phenotypic equilibrium as probabilistic convergence in
  multi-phenotype cell population dynamics.
\newblock {\em PLOS ONE 12}, 2 (2017), e0170916.

\bibitem{kaimanovich2002boundary}
{\sc Kaimanovich, V.~A., and Woess, W.}
\newblock Boundary and entropy of space homogeneous markov chains.
\newblock {\em Annals of probability\/} (2002), 323--363.

\bibitem{kang2014flexibility}
{\sc Kang, Y., Gu, C., Yuan, L., Wang, Y., Zhu, Y., Li, X., Luo, Q., Xiao, J.,
  Jiang, D., Qian, M., et~al.}
\newblock Flexibility and symmetry of prokaryotic genome rearrangement reveal
  lineage-associated core-gene-defined genome organizational frameworks.
\newblock {\em mBio 5}, 6 (2014), e01867--14.

\bibitem{ma2015universal}
{\sc Ma, Y.-A., and Qian, H.}
\newblock Universal ideal behavior and macroscopic work relation of linear
  irreversible stochastic thermodynamics.
\newblock {\em New Journal of Physics 17}, 6 (2015), 065013.

\bibitem{maskin2014arrow}
{\sc Maskin, E., and Sen, A.}
\newblock {\em The Arrow impossibility theorem}.
\newblock Columbia University Press, 2014.

\bibitem{niu2015phenotypic}
{\sc Niu, Y., Wang, Y., and Zhou, D.}
\newblock The phenotypic equilibrium of cancer cells: {From} average-level
  stability to path-wise convergence.
\newblock {\em Journal of Theoretical Biology 386\/} (2015), 7--17.

\bibitem{qian2020counting}
{\sc Qian, H., and Cheng, Y.-C.}
\newblock Counting single cells and computing their heterogeneity: from
  phenotypic frequencies to mean value of a quantitative biomarker.
\newblock {\em Quantitative Biology 8}, 2 (2020), 172--176.

\bibitem{qian2019ternary}
{\sc Qian, H., Cheng, Y.-C., and Thompson, L.~F.}
\newblock Ternary representation of stochastic change and the origin of entropy
  and its fluctuations.
\newblock {\em arXiv preprint arXiv:1902.09536\/} (2019).

\bibitem{qian2020kinematic}
{\sc Qian, H., Cheng, Y.-C., and Yang, Y.-J.}
\newblock Kinematic basis of emergent energetics of complex dynamics.
\newblock {\em EPL (Europhysics Letters) 131}, 5 (2020), 50002.

\bibitem{wang2018some}
{\sc Wang, Y.}
\newblock {\em Some Problems in Stochastic Dynamics and Statistical Analysis of
  Single-Cell Biology of Cancer}.
\newblock {Ph.D. thesis}, University of Washington, 2018.

\bibitem{wang2022impossibility}
{\sc Wang, Y.}
\newblock Impossibility results about inheritance and order of death.
\newblock {\em PLOS ONE 17}, 11 (2022), e0277430.

\bibitem{wang2022two}
{\sc Wang, Y.}
\newblock Two metrics on rooted unordered trees with labels.
\newblock {\em Algorithms for Molecular Biology 17}, 1 (2022), 1--17.

\bibitem{wang2023longest}
{\sc Wang, Y.}
\newblock Longest common subsequence algorithms and applications in determining
  transposable genes.
\newblock {\em arXiv preprint arXiv:2301.03827\/} (2023).

\bibitem{wang2022modelling}
{\sc Wang, Y., Dessalles, R., and Chou, T.}
\newblock Modelling the impact of birth control policies on {China}'s
  population and age: effects of delayed births and minimum birth age
  constraints.
\newblock {\em Royal Society Open Science 9}, 6 (2022), 211619.

\bibitem{wang2022ar}
{\sc Wang, Y., and He, S.}
\newblock Inference on autoregulation in gene expression.
\newblock {\em arXiv preprint arXiv:2201.03164\/} (2022).

\bibitem{wang2020biological}
{\sc Wang, Y., Kropp, J., and Morozova, N.}
\newblock Biological notion of positional information/value in morphogenesis
  theory.
\newblock {\em International Journal of Developmental Biology 64}, 10-11-12
  (2020), 453--463.

\bibitem{wang2020model}
{\sc Wang, Y., Minarsky, A., Penner, R., Soul{\'e}, C., and Morozova, N.}
\newblock Model of morphogenesis.
\newblock {\em Journal of Computational Biology 27}, 9 (2020), 1373--1383.

\bibitem{wang2022discrete}
{\sc Wang, Y., Mistry, B.~A., and Chou, T.}
\newblock Discrete stochastic models of selex: Aptamer capture probabilities
  and protocol optimization.
\newblock {\em The Journal of Chemical Physics 156}, 24 (2022), 244103.

\bibitem{wang2020mathematical}
{\sc Wang, Y., and Qian, H.}
\newblock Mathematical representation of {Clausius'} and {Kelvin's} statements
  of the second law and irreversibility.
\newblock {\em Journal of Statistical Physics 179}, 3 (2020), 808--837.

\bibitem{wang2020causal}
{\sc Wang, Y., and Wang, L.}
\newblock Causal inference in degenerate systems: {An} impossibility result.
\newblock In {\em International Conference on Artificial Intelligence and
  Statistics\/} (2020), PMLR, pp.~3383--3392.

\bibitem{wang2022inference}
{\sc Wang, Y., and Wang, Z.}
\newblock Inference on the structure of gene regulatory networks.
\newblock {\em Journal of Theoretical Biology 539\/} (2022), 111055.

\bibitem{wang2021inference}
{\sc Wang, Y., Zhang, B., Kropp, J., and Morozova, N.}
\newblock Inference on tissue transplantation experiments.
\newblock {\em Journal of Theoretical Biology 520\/} (2021), 110645.

\bibitem{wang2021measuring}
{\sc Wang, Y., and Zheng, Z.}
\newblock Measuring policy performance in online pricing with offline data.
\newblock {\em Available at SSRN 3729003\/} (2021).

\bibitem{wang2023online}
{\sc Wang, Y., Zheng, Z., and Shen, Z.-J.~M.}
\newblock Online pricing with polluted offline data.
\newblock {\em Available at SSRN 4320324\/} (2023).

\bibitem{wang2023multiple}
{\sc Wang, Y., Zhou, J.~X., Pedrini, E., Rubin, I., Khalil, M., Qian, H., and
  Huang, S.}
\newblock {Multiple phenotypes in HL60 leukemia cell population}.
\newblock {\em arXiv preprint arXiv:2301.03782\/} (2023).

\bibitem{wolpert2016positional}
{\sc Wolpert, L.}
\newblock Positional information and pattern formation.
\newblock {\em Current topics in developmental biology 117\/} (2016), 597--608.

\bibitem{xia2023age}
{\sc Xia, M., Li, X., and Chou, T.}
\newblock An age-structured {Lotka-Volterra} model and the emergence of
  overcompensation.
\newblock {\em Bulletin of the American Physical Society\/} (2023).

\bibitem{yang2021potentials}
{\sc Yang, Y.-J., and Cheng, Y.-C.}
\newblock Potentials of continuous {Markov} processes and random perturbations.
\newblock {\em Journal of Physics A: Mathematical and Theoretical 54}, 19
  (2021), 195001.

\bibitem{ye2016stochastic}
{\sc Ye, F. X.-F., Wang, Y., and Qian, H.}
\newblock Stochastic dynamics: {Markov} chains and random transformations.
\newblock {\em Discrete \& Continuous Dynamical Systems-B 21}, 7 (2016), 2337.

\bibitem{zhao2016part}
{\sc Zhao, J., Zhou, Y., Zhang, X., and Chen, L.}
\newblock Part mutual information for quantifying direct associations in
  networks.
\newblock {\em Proceedings of the National Academy of Sciences 113}, 18 (2016),
  5130--5135.

\bibitem{zhou2014multi}
{\sc Zhou, D., Wang, Y., and Wu, B.}
\newblock A multi-phenotypic cancer model with cell plasticity.
\newblock {\em Journal of Theoretical Biology 357\/} (2014), 35--45.

\bibitem{zhou2021dissecting}
{\sc Zhou, P., Wang, S., Li, T., and Nie, Q.}
\newblock Dissecting transition cells from single-cell transcriptome data
  through multiscale stochastic dynamics.
\newblock {\em Nature communications 12}, 1 (2021), 1--15.

\end{thebibliography}

\end{document}